\documentclass[final]{siamltex}

\usepackage[parfill]{parskip}    
\usepackage{amsmath,amssymb,latexsym,enumerate,mathrsfs}
\usepackage{hyperref}
\usepackage{array}
\usepackage{times}
\usepackage[]{mdframed}
\usepackage{epstopdf}
\usepackage{subfig}
\usepackage{graphicx}
\usepackage{hyperref} 
\usepackage{multicol}
\usepackage{framed}
\usepackage{moreverb}
\usepackage{marvosym}
\usepackage{color,soul}
\definecolor{lightblue}{rgb}{.90,.95,1}
\sethlcolor{yellow}

\DeclareMathOperator*{\argmin}{arg\,min}

\newtheorem{remark}{Remark}

\title{Further extensions on the successive  approximation method for hierarchical optimal control problems and its application to learning}

\author{Getachew K. Befekadu}

\begin{document}
\maketitle

\renewcommand{\thefootnote}{\arabic{footnote}}

\begin{abstract}
In this paper, further extensions of the result of the paper ``{\bf A successive approximation method in functional spaces for hierarchical optimal control problems and its application to learning, arXiv:2410.20617 [math.OC], 2024}'' concerning a class of learning problem of point estimations for modeling of high-dimensional nonlinear functions are given. In particular, we present two viable extensions within the nested algorithm of the successive approximation method for the hierarchical optimal control problem, that provide better convergence property and computationally efficiency, which ultimately leading to an optimal parameter estimate. The first extension is mainly concerned with the convergence property of the steps involving how the two agents, i.e., the ``{\it leader}'' and the ``{\it follower},'' update their admissible control strategies, where we introduce augmented Hamiltonians for both agents and we further reformulate the admissible control updating steps as as sub-problems within the nested algorithm of the hierarchical optimal control problem that essentially provide better convergence property. Whereas the second extension is concerned with the computationally efficiency of the steps involving how the agents update their admissible control strategies, where we introduce intermediate state variable for each agent and we further embed the intermediate states within the optimal control problems of the ``{\it leader}'' and the ``{\it follower},'' respectively, that further lend the admissible control updating steps to be fully efficient time-parallelized within the nested algorithm of the hierarchical optimal control problem.
\end{abstract}

\begin{keywords} 
Augmented Hamiltonian, controllability, generalization, hierarchical optimal control problems, learning problem, modeling of nonlinear functions, modified-successive approximation method, Pontryagin's maximum principle, reachability, regularization, Stackelberg's optimization, time-parallelized distributed computation.
\end{keywords}

\section{Introduction} \label{S1}
Recently, an interesting connection between a class of learning problem of point estimations for modeling of high-dimensional nonlinear functions and that of a hierarchical optimal control problem has been addressed in the paper \cite{r1}. In fact, the paper has established an evidential connection between the optimal control problem and an inverse problem whose optimality conditions are the direct consequence of the {\it Pontryagin's maximum principle}. Moreover, the paper also presented a nested algorithm, based on successive approximation methods, for solving numerically the corresponding hierarchical optimal control problem that ultimately leading to an optimal parameter estimate. In this paper, we present two fairly practical extensions within the nested algorithm of the successive approximation method for the hierarchical optimal control problem that provide better convergence and computationally efficiency. In particular, the first extension is mainly concerned with the convergence property of the steps involving how the two agents, i.e., the ``{\it leader}'' and the ``{\it follower},'' update their admissible control strategies, where we specifically introduce augmented Hamiltonians for both agents and we further reformulate the admissible control updating steps as as sub-problems within the nested algorithm of the hierarchical optimal control problem that essentially provide better convergence property. Whereas the second extension is concerned with the computationally efficiency of the steps involving how the agents update their admissible control strategies, where we specifically introduce intermediate state variable for each agent and we further embed the corresponding intermediate state variables within the optimal control problems of the ``{\it leader}'' and the ``{\it follower},'' respectively, that further lend the admissible control updating steps to be fully efficient time-parallelized within the nested algorithm of the hierarchical optimal control problem. In this paper, our intent is to provide a more viable computational frameworks, rather than considering any specific numerical problem. Some numerical works using the proposed computational frameworks have been done and detailed results will be presented elsewhere.\footnote{Here we emphasize that we are not proposing a new optimization framework, but rather -- based on the evidential connection between the optimal control theory and inverse problems -- we present extensions within the nested algorithm of the optimal control problem that provide better convergence and computationally efficiency.}

The remainder of this paper is organized as follows. In Section~\ref{S2}, we restate some of the core concepts and assumptions that have been presented in the paper \cite{r1}. In Section~\ref{S3}, we present our main results, where we provide two further extensions within the nested algorithm of the successive approximation method for the hierarchical optimal control problem, that provide better convergence property and computationally efficiency. For the sake of completeness, we also provide the original nested algorithm (which is presented in the paper \cite{r1}) in the Appendix section.

\section{Problem statement} \label{S2}
In this section, we restate some of the core concepts and assumptions that have been presented in the paper \cite{r1}. In fact, the paper presented a mathematical construct how to account appropriately for both generalization and regularization based on an evidential connection between a class of learning problem and that of a hierarchical optimal control problem. 

Here, the problem statement consists of the following core concepts and general assumptions:

\begin{enumerate} [(i).]
\item {\it Datasets}: We are given two datasets, i.e., $\mathcal{Z}^{(k)} = \bigl\{ (x_i^{(k)}, y_i^{(k)})\bigr\}_{i=1}^{m_k}$, each with data size of $m_k$, for $k=1, 2$. The datasets $\bigl\{ \mathcal{Z}^{(1)} \bigr\}_{k=1}^{2}$ may be generated from a given original dataset $\mathcal{Z}^{(0)} =\bigl\{ (x_i^{(0)}, y_i^{(0)})\bigr\}_{i=1}^{m_0}$ by means of bootstrapping with/without replacement. Here, we assume that the first dataset $\mathcal{Z}^{(1)}=\bigl\{ (x_i^{(1)}, y_i^{(1)})\bigr\}_{i=1}^{m_1}$ will be used for model training purpose, while the second dataset $\mathcal{Z}^{(2)} = \bigl\{ (x_i^{(2)}, y_i^{(2)})\bigr\}_{i=1}^{m_2}$ will be used for evaluating the quality of the estimated model parameter.
\item {\it Learning via gradient systems}: We are tasked to find for a parameter $\theta \in \Theta$, from a finite-dimensional parameter space $\mathbb{R}^p$ (i.e., $\Theta \subset \mathbb{R}^p$), such that the function $h_{\theta}(x) \in \mathcal{H}$, i.e., from a given class of hypothesis function space $\mathcal{H}$, describes best the corresponding model training dataset as well as predicts well with reasonable expectation on a different model validation dataset. Here, the search for an optimal parameter $\theta^{\ast} \in \Gamma \subset \mathbb{R}^p$ can be associated with the {\it steady-state solution} to the following gradient system, whose {\it time-evolution} is guided by the model training dataset $\mathcal{Z}^{(1)}$, i.e.,
\begin{align}
 \dot{\theta}(t) = - \nabla J_0(\theta(t),\mathcal{Z}^{(1)}), \quad \theta(0) = \theta_0, \label{Eq2.1}
\end{align}
with $J_0(\theta, \mathcal{Z}^{(1)}) = \frac{1}{m_1} \sum\nolimits_{i=1}^{m_1} {\ell} \bigl(h_{\theta}(x_i^{(1)}), y_i^{(1)} \bigr)$, where $\ell$ is a suitable loss function that quantifies the lack-of-fit between the model and the datasets.
\item {\it Optimal control problem}: Instead of using directly the gradient system of Equation~\eqref{Eq2.1}, we consider the following controlled-gradient system
\begin{align}
 \dot{\theta}^{u}(t) = - \nabla J_0(\theta^{u}(t),\mathcal{Z}^{(1)}) + u(t), \quad \theta^{u}(0) = \theta_0, \label{Eq2.2}
 \end{align}
where $u(t)$ is a control strategy to be determined from a class of admissible control space $\mathcal{U}$, i.e., $\mathcal{U} = L^{\infty}([0,T], U)$ is the space of measurable and essentially bounded functions from $[0,T]$ to $U \subset \mathbb{R}^p$.\footnote{A function $u \colon [0,T] \to U \subset \mathbb{R}^p$ is said to be essentially bounded, if there exists a set $\omega \subset [0,T]$ of measure zero, such that $u$ is bounded on $[0,T] \setminus \omega.$} 
 Then, we search for an optimal control strategy $u^{\ast}(t) \in \mathcal{U}$, so that we would like to achieve the following two objectives:
 \begin{enumerate} [(1).]
\item {\it Controllability-type objective}: Suppose that we are given a target set $\Gamma \in \Theta$ that may depend on the model validation dataset $\mathcal{Z}^{(2)}$, where such a target set can be used as a criterion for evaluating the quality of the estimated model parameter. Then, the final estimated parameter, at some fixed time $T$, is expected to reach this target set $\Gamma$, i.e.,
\begin{align}
 \theta^{u^{\ast}}(T) \in \Gamma, \label{Eq2.3}
 \end{align}
starting from an initial point $\theta^{u}(0) = \theta_0$.
\item {\it Regularization-type objective }: Under certain conditions (see below the general assumptions), we would like to ensure the estimated parameter trajectory $\theta^{u}(t) \in \Theta^{M}$, with $u(t) \in \mathcal{U}$, where $\Theta^{M}$ is a $p$-dimensional manifold, to satisfy some regularization property over a finite time interval $[0, T]$. 
\end{enumerate}
\item {\it General assumptions:} We assume the following conditions\footnote{Note that these conditions are sufficient for the existence of a nonempty compact reachable set $\mathcal{R}(\theta_0)$, for some admissible controls on $[0,T]$ that belongs to $\mathcal{U}$, starting from an initial point $\theta^{u}(0)=\theta_0$ (e.g., see \cite{r2} for related discussions on the Filippov's theorem providing a sufficient condition for compactness of the reachable set). Moreover, the controllability-type problem is solvable, if there exists $u^{\ast}(t) \in \mathcal{U}$ such that $\mathcal{R}(\theta^{u^{\ast}}(T)) \cap \Gamma \neq \emptyset$.}
 \begin{enumerate} [(a).]
 \item The set $U \subset \mathbb{R}^{p}$ is compact and the final time $T$ is fixed.
 \item The set $\tilde{\Theta} = \left \{- \nabla J_0(\theta,\mathcal{Z}^{(1)}) + u \,\bigl\vert\, u \in U \right\}$ is convex for every $\theta \in \Theta^{M}$.
 \item For every $u(t) \in \mathcal{U}$, the solution of the controlled-gradient system $\dot{\theta}^{u}(t) = - \nabla J_0(\theta^{u}(t),\mathcal{Z}^{(1)}) + u(t)$, with $\theta^{u}(0) = \theta_0$, is defined on the whole interval $[0,T]$ and belongs to $\Theta^{M}$.
 \end{enumerate}
\end{enumerate}
In order to make the above optimal control problem mathematically more precise, we consider the following hierarchical cost functionals:
\begin{align}
J_1[u] &= \int_0^T \frac{1}{2} \bigl \Vert \theta^{u}(t) \bigr \Vert^2 dt \notag\\
             & \text{s.t.} \quad \theta^{u}(T) \in \Gamma \label{Eq2.4}
\end{align}
and
\begin{align}
J_2[u] = \int_0^T \left \{ \frac{\alpha}{2} \bigl\Vert \theta^{u}(t) \bigr \Vert^2 + \frac{\beta}{2} \bigl\Vert u(t) \bigr\Vert^2 \right\} dt, \quad \alpha > 0, \quad \beta > 0. \label{Eq2.5}
\end{align}
In general, finding an optimal control strategy $u^{\ast}(t) \in \mathcal{U}$, i.e., a {\it Pareto-optimal solution}, that simultaneously minimizes the above two cost functionals is not an easy problem. However, using the notion of {\it Stackelberg's optimization} (see \cite{r3}, \cite{r4} or \cite{r5} for additional discussions), we provide conditions on the existence of admissible optimal control strategies for such a hierarchical optimal control problem. Here, the main idea is to partition the admissible control strategy space $\mathcal{U}$ into two open subspaces $\mathcal{U}_1$ and $\mathcal{U}_2$, with smooth boundaries ({\it up to a set of measurable} $\mathcal{U}$, with $\mathcal{U}_1 \cap \mathcal{U}_2 = \varnothing$), that are compatible with two abstract agents, namely., a ``{\it leader}'' (which is responsible for the controllability-type problem) and that of a ``{\it follower}'' (which is associated with the regularization-type problem). 

Then, the above optimal control problem will become searching for a pair of admissible optimal control strategies $(u_1^{\ast}(t),u_2^{\ast}(t)) \in \mathcal{U}_1 \times \mathcal{U}_2$ that steers the controlled-gradient system
\begin{align}
 \dot{\theta}^{u}(t) = - \nabla J_0(\theta^{u}(t),\mathcal{Z}^{(1)}) + u_1(t) \chi_{\mathcal{U}_1} + u_2(t) \chi_{\mathcal{U}_2}, \quad \theta^{u}(0) = \theta_0 \label{Eq2.6}
\end{align}
so as to achieve the overall objectives of Equations~\eqref{Eq2.4} and \eqref{Eq2.5} that ultimately leading to an optimal parameter estimate $\theta^{u^{\ast}}(T) \in \Gamma$.\footnote{Here, we use the notation $\chi_{\mathcal{U}_i}$, for $i=1,2$, to denote the characteristic function for $\mathcal{U}_i$, where the admissible control $u_i(t)$ is the restriction of the admissible control $u$ to $\mathcal{U}_i$.} 

In what follows, we provide the optimality conditions for both agents, i.e., the ``follower'' and  the ``leader,'' where such conditions are the direct consequence of the {\it Pontryagin's maximum principle} (see \cite{r1} for detailed proofs).

\subsection*{Optimality statement for the follower}
Suppose that an admissible control strategy $\tilde{u}_1(t) \in \mathcal{U}_1$ for the ``{\it leader}'' is given. Then, the solution for the optimal control problem w.r.t. the ``{\it follower,}'' i.e.,
\begin{align}
J_2[u_2] &= \int_0^T \left \{ \frac{\alpha}{2} \bigl \Vert \theta^{u_2}(t) \bigr \Vert^2  + \frac{\beta}{2} \bigl \Vert u_2(t) \chi_{\mathcal{U}_2} \bigr \Vert^2 \right\} dt \quad \to \quad \min_{u_2() \in \mathcal{U}_2} \label{Eq2.7}\\
             & \text{s.t.} \notag \\
             \quad \dot{\theta}^{u_2}(t) &= -\nabla J_0(\theta^{u_2}(t), \mathcal{Z}^{(1)}) + \tilde{u}_1(t) \chi_{\mathcal{U}_1} + u_2(t) \chi_{\mathcal{U}_2}, \quad \theta^{u_2}(0) = \theta_0, \label{Eq2.8}
\end{align}
satisfies the following Euler-Lagrange critical point equations:
\begin{enumerate} [(i).]
\item the forward-equation, i.e., the state estimation dynamics\footnote{Notice that the solution for the state estimation trajectory $\theta^{u_2}(\cdot) \in \Theta^{M}$ and it depends on both the strategies of the ``{\it follower}'' and the ``{\it leader}.''}
\begin{align}
 \dot{\theta}^{u_2}(t) &= \frac{\partial H_2(\theta^{u_2}, p_2, u_2)}{\partial p_2} \notag \\
                                  &=-\nabla J_0(\theta^{u_2}(t), \mathcal{Z}^{(1)}) + \tilde{u}_1(t) \chi_{\mathcal{U}_1} + u_2(t) \chi_{\mathcal{U}_2}, \quad \theta^{u_2}(0) = \theta_0, \label{Eq2.9}
\end{align}
\item the backward-equation, i.e., the adjoint state equation w.r.t. the ``{\it follower,}'' 
\begin{align}
  \dot{p}_2(t) = - \frac{\partial H_2(\theta^{u_2}, p_2, u_2)}{\partial \theta^{u_2}}, \quad p_2(T) = 0, \label{Eq2.10}
\end{align}
\item the extremum condition
\begin{align}
  \frac{\partial H_2(\theta^{u_2}, p_2, u_2)}{\partial u_2} = 0 \quad {\text on} \quad [0,T], \label{Eq2.11}
\end{align}
\end{enumerate}
where $H_2$ is the Hamiltonian equation w.r.t. the ``{\it follower}'' and it is given by
\begin{align}
H_2(\theta^{u_2}, p_2, u_2) &= \bigl \langle -\nabla J_0(\theta^{u_2}(t), \mathcal{Z}^{(1)}) + \tilde{u}_1(t) \chi_{\mathcal{U}_1} + u_2(t) \chi_{\mathcal{U}_2}, p_2 \bigr\rangle  \notag \\
& \quad \quad + \frac{\alpha}{2} \bigl\Vert \theta^{u_2}(t) \bigr\Vert^2  + \frac{\beta}{2} \bigl\Vert u_2(t) \chi_{\mathcal{U}_2} \bigr \Vert^2. \label{Eq2.12}
\end{align}

\begin{remark} \label{R1}
Note that the above Euler-Lagrange critical point equations necessitate the following argument, where the ``{\it follower}'' is required to respond optimally to the strategy of the ``{\it leader},'' i.e., the optimal strategy $u_2 \in \mathcal{U}_2$ for the ``{\it follower}'' is a functional on $[0, T]$ that depends on $u_1(\cdot) \in \mathcal{U}_1$. Hence, such a correspondence problem between the ``{\it leader}'' and that of the ``{\it follower}'' admits a unique functional mapping $\mathcal{F} [~]$, i.e.,
\begin{align}
u_2(t) = \mathcal{F}[u_1(t) \chi_{\mathcal{U}_1}] \in \mathcal{U}_2 \quad \left (\text{i.e.,} \quad \mathcal{F} \colon \mathcal{U}_1 \mapsto \mathcal{U}_2 \quad \text{on}\quad [0,T] \right). \label{Eq2.13}
\end{align}
 \end{remark}
 
\subsection*{Optimality statement for the leader} Suppose that the critical point conditions in Equations~\eqref{Eq2.9}, \eqref{Eq2.10} and \eqref{Eq2.11}hold true. Then, the solution for the optimal control problem w.r.t. the ``{\it leader,}'' i.e.,
\begin{align}
J_1[u_1] &= \int_0^T \frac{1}{2} \bigl \Vert \theta^{u_1}(t) \bigr \Vert^2 dt \quad \to \quad \min_{u_1() \in \mathcal{U}_1} \label{Eq2.14}\\
             & \text{s.t.} \quad\quad \Phi(\theta^{u_1}(T),\mathcal{Z}^{(2)}) = z, \label{Eq2.15} \\
             \quad \dot{\theta}^{u_1}(t) &= -\nabla J_0(\theta^{u_1}(t), \mathcal{Z}^{(1)}) + u_1(t) \chi_{\mathcal{U}_1} + \mathcal{F}[u_1(t)\chi_{\mathcal{U}_1}] \chi_{\mathcal{U}_2}, \quad \theta^{u_1}(0) = \theta_0, \label{Eq2.16} 
\end{align}
satisfies the following Euler-Lagrange critical point equations:
\begin{enumerate} [(i).]
\item the forward-equation, i.e., the state estimation dynamics
\begin{align}
 \dot{\theta}^{u_1}(t) &= \frac{\partial H_1(\theta^{u_1}, p_1, u_1)}{\partial p_1} \notag \\
                                  &=-\nabla J_0(\theta^{u_1}(t), \mathcal{Z}^{(1)}) + u_1(t) \chi_{\mathcal{U}_1} + \mathcal{F}[u_1(t)\chi_{\mathcal{U}_1}] \chi_{\mathcal{U}_2}, \quad \theta^{u_1}(0) = \theta_0, \label{Eq2.17}
\end{align}
\item the backward-equation, i.e., the adjoint state equation w.r.t. the ``{\it leader},'' 
\begin{align}
 \dot{p}_1(t) = - \frac{\partial H_1(\theta^{u_1}, p_1, u_1)}{\partial \theta^{u_1}}, \quad p_1(T) = -\frac{\partial \Phi(\theta,\mathcal{Z}^{(1)})} {\partial \theta} \biggr \vert_{\theta=\theta^{u_1}(T)}, \label{Eq2.18}
                                  \end{align}
\item the extremum condition
\begin{align}
 \frac{\partial H_1(\theta^{u_1}, p_1, u_1)}{\partial u_1} = 0 \quad {\text on} \quad [0,T],  \label{Eq2.19}
\end{align}
\end{enumerate}
where $H_1$ is the Hamiltonian equation w.r.t. the ``{\it leader}'' and it is given by
\begin{align}
H_1(\theta^{u_1}, p_1, u_1) &= \bigl \langle -\nabla J_0(\theta^{u_1}(t), \mathcal{Z}^{(1)}) + u_1(t) \chi_{\mathcal{U}_1} + \mathcal{F}[u_1(t)\chi_{\mathcal{U}_1} ] \chi_{\mathcal{U}_2}, p_1 \bigr\rangle  \notag \\
& \quad \quad + \frac{1}{2} \bigl \Vert \theta^{u_1}(t) \bigr \Vert^2. \label{Eq2.20}
\end{align}

\section{Main results} \label{S3}
In this section, we present our main results, where we present two viable extensions -- based on the modified successive approximation methods and the intermediate state methods -- within the nested algorithm of the successive approximation method for the hierarchical optimal control problem.

\subsection{Modified successive approximation methods} \label{S3.1}
In what follows, we present a modified successive approximation method within the nested algorithm for the solution of the hierarchical optimal control problem. Here, the modification involves two phases within the nested algorithm: (i) we first introduce an augmented Hamiltonian w.r.t. the ``{\it follower}'' and we proceed by modifying how the ``{\it follower}'' updates its admissible control strategy, which involves minimizing the corresponding augmented Hamiltonian. Note that this part of the modification only affects ${\rm Step~3}$ and ${\rm Step~4}$ in the original nested algorithm (see also {\rm ALGORITHM - O} in the Appendix section). (ii) We similarly introduce another augmented Hamiltonian w.r.t. the ``{\it leader}'' and we then proceed by modifying how the ``{\it leader}'' updates its admissible control strategy, which also involves minimizing the corresponding augmented Hamiltonian. Note that the modification at this time only affects ${\rm Step~6}$ and ${\rm Step~7}$ in the original nested algorithm.

\subsubsection{Modification w.r.t. the follower's admissible control updating} \label{S3.1.1}
Let us assume that we have chosen a nominal control strategy $\bar{u}_2(t) \in \mathcal{U}_2$ w.r.t. the ``{\it follower}'' and we further obtained the corresponding nominal trajectory $\theta^{\bar{u}_2}(t)$ by solving the system dynamics in the forward direction. Next, we introduce an augmented Hamiltonian $\tilde{H}_2$ w.r.t. the ``{\it follower}'' as follows, with $\gamma_2 \in [0, 1)$,
\begin{align}
\tilde{H}_2(\theta^{u_2}, p_2, \bar{u}_2, u_2) &= H_2(\theta^{u_2}, p_2, u_2) + \frac{\gamma_2}{2}\biggl \Vert \frac{\partial H_2(\theta^{u_2}, p_2, u_2)}{\partial p_2} - \frac{\partial H_2(\theta^{\bar{u}_2}, p_2, \bar{u}_2)}{\partial p_2}\biggr \Vert^2 \notag\\
& \quad \quad + \frac{\gamma_2}{2}\biggl \Vert \frac{\partial H_2(\theta^{u_2}, p_2, u_2)}{\partial \theta^{u_2}} - \frac{\partial H_2(\theta^{\bar{u}_2}, p_2, \bar{u}_2)}{\partial \theta^{\bar{u}_2}}\biggr \Vert^2 \label{Eq3.1}
\end{align}
Then, based on the newly defined augmented Hamiltonian, the updating step for the control strategy of the ``{follower}'' will involve solving an auxiliary minimization problem, i.e., 
\begin{align}
u_2^{\ast} \in \argmin_{u_2 \in \mathcal{U}_2} \tilde{H}_2(\theta^{u_2}, p_2, \bar{u}_2, u_2), \quad \forall_{t \in [0, T]}, \label{Eq3.2}
\end{align}
where the improvement in the control strategy $\delta u_2$ will involve solving the first and second order variations of $J_2$ in the auxiliary minimization problem (e.g., see \cite{r7} for additional discussions). Note that part of the modification w.r.t. the ``{follower}''  in the nested algorithm will necessitate the following steps:

\begin{itemize}
\item[{ i.}] Using the admissible control pairs $\bigl(u_1^{(n-1)}(t), u_2^{(n-1)}(t)\bigr)$, solve the forward and backward-equations w.r.t. the system dynamics of the ``{follower}''.
\item[{ ii.}] Update the control strategy of the ``{follower}'' by solving the following auxiliary minimization problem
\begin{align}
u_2^{(n)} \in \argmin_{u_2 \in \mathcal{U}_2} \tilde{H}_2(\theta^{u_2}, p_2, u_2^{(n-1)}, u_2), \quad \forall_{t \in [0, T]}. \label{Eq3.3}
\end{align}
\end{itemize}
 
\subsubsection{Modification w.r.t. the leader's admissible control updating} \label{S3.1.2} Here, we similarly introduce an augmented Hamiltonian $\tilde{H}_1$ w.r.t. the ``{\it leader}'' as follows, with $\gamma_1 \in [0, 1)$,
\begin{align}
\tilde{H}_1(\theta^{u_1}, p_1, \bar{u}_1, u_1) &= H_1(\theta^{u_1}, p_1, u_1) + \frac{\gamma_1}{2}\biggl \Vert \frac{\partial H_1(\theta^{u_1}, p_1, u_1)}{\partial p_1} - \frac{\partial H_1(\theta^{\bar{u}_1}, p_1, \bar{u}_1)}{\partial p_1}\biggr \Vert^2 \notag\\
& \quad \quad + \frac{\gamma_1}{2}\biggl \Vert \frac{\partial H_1(\theta^{u_1}, p_1, u_1)}{\partial \theta^{u_1}} - \frac{\partial H_1(\theta^{\bar{u}_1}, p_1, \bar{u}_1)}{\partial \theta^{\bar{u}_1}}\biggr \Vert^2 \label{Eq3.4}
\end{align}
Then, using the newly defined augmented Hamiltonian, the updating step for the control strategy of the ``{leader}'' will involve solving a separate auxiliary minimization problem, i.e., 
\begin{align}
u_1^{\ast} \in \argmin_{u_1 \in \mathcal{U}_1} \tilde{H}_1(\theta^{u_1}, p_1, \bar{u}_1, u_1), \quad \forall_{t \in [0, T]}, \label{Eq3.5}
\end{align}
where the improvement in the control strategy $\delta u_1$ will involve solving the first and second order variations of $J_1$ in the auxiliary minimization problem. Note that part of the modification w.r.t. the ``{leader}'' in the nested algorithm will necessitate the following steps:

\begin{itemize}
\item[{ i.}] Using the admissible control pairs $\bigl(u_1^{(n-1)}(t), u_2^{(n)}(t)\bigr)$, solve the forward and backward-equations w.r.t. the system dynamics of the ``{leader}''.
\item[{ ii.}] Update the control strategy of the ``{leader}'' by solving the following auxiliary minimization problem
 \begin{align}
u_1^{(n)} \in \argmin_{u_1 \in \mathcal{U}_1} \tilde{H}_1(\theta^{u_1}, p_1, u_1^{(n-1)}, u_1), \quad \forall_{t \in [0, T]}. \label{Eq3.6}
\end{align}
\end{itemize}

\begin{remark} \label{R2}
Here, we remark that, when both $\gamma_1 = 0$ and $\gamma_2 = 0$, the augmented Hamiltonians will be the same as the original Hamiltonians w.r.t. the ``{follower}'' and the ``{leader}'' and method reduces to the original nested algorithm based on the successive approximation method for the solution of the hierarchical optimal control problem.
\end{remark}

\subsubsection*{Algorithm -- Modified successive approximation methods}
Note that if we update the original nested algorithm based on the above two modifications w.r.t. the ``{follower}'' and that of the ``{leader}''. Then, we have the following modified nested algorithm based on the modified successive approximation methods that provides better convergence property.

{\rm \footnotesize

{\bf ALGORITHM - 1:} Modified Nested Algorithm Based on the Modified Successive Approximation Methods
\begin{itemize}
\item[{\bf 0.}] {\bf Initialize:} Start with any admissible control strategy $u_1^{(n-1)}(t) \in \mathcal{U}_1$ (i.e., an initial guess control function on $[0,T]$) for the ``{\it leader}.'' 
\item[{\bf 1.}] {\bf The regularization-type problem:} Choose an admissible control strategy $u_2^{(n-1)}(t) \in \mathcal{U}_2$ (i.e., an initial guess control function on $[0,T]$) for the ``{\it follower}.'' 
\item[{\bf 2.}] Using the admissible control pairs $\bigl(u_1^{(n-1)}(t), u_2^{(n-1)}(t)\bigr)$, solve the forward and backward-equations w.r.t. the system dynamics of the ``{\it follower},'' i.e., Equations~\eqref{Eq2.9} and \eqref{Eq2.10}.
 \item[{\bf 3.}] Update the new admissible control for the ``{\it follower}, i.e., $u_2^{(n)}(t)$, using
\begin{align*}
u_2^{(n)} \in \argmin_{u_2 \in \mathcal{U}_2} \tilde{H}_2(\theta^{u_2}, p_2, u_2^{(n-1)}, u_2), \quad \forall_{t \in [0, T]}
\end{align*}
 \item[{\bf 4.}] {\bf The controllability-type problem:} Using the admissible control strategy pairs $\bigl(u_1^{(n-1)}(t), u_2^{(n)}(t)\bigr)$, with the updated control strategy $u_2^{(n)}(t)$, solve the forward and backward-equations w.r.t. the system dynamics of the ``{\it leader,}'' i.e., Equations~\eqref{Eq2.17} and \eqref{Eq2.18}.
 \item[{\bf 5.}] Update the new admissible control for the ``{\it leader},''  i.e., $u_1^{(n)}(t)$, using
\begin{align*}
u_1^{(n)} \in \argmin_{u_1 \in \mathcal{U}_1} \tilde{H}_1(\theta^{u_1}, p_1, u_1^{(n-1)}, u_1), \quad \forall_{t \in [0, T]}
\end{align*}
\item[{\bf 6.}] With the updated admissible control strategy pairs $\bigl(u_1^{(n)}(t), u_2^{(n)}(t)\bigr)$, repeat Steps $2$ through $5$, until convergence, i.e., $\bigl\Vert \partial \tilde{H}_1/ \partial u_1 \bigr\Vert \le \epsilon_{\rm tol}$, for some error tolerance $\epsilon_{\rm tol} > 0$.
 \item[{\bf 7.}] {\bf Output:} Return the optimal estimated parameter value $\theta^{\ast} = \theta^{u^{(n)}}(T)$, with $u^{(n)}(t) = (u_2^{(n)}(t), u_2^{(n)}(t))$.
\end{itemize}}

\subsection{Intermediate state methods} \label{S3.2}
In what follows, we present the extension based on the intermediate state method, i.e., a time-parallelized method, into two phases. Here, we assume that the solutions for both the forward and backward equations corresponding to the ``{follower}'' and the ``{leader}'' are first obtained on the time interval $[0,T]$, i.e., we have solved the corresponding solutions at $t_0=0$, $t_1=\delta$, $t_2=2\delta$, \ldots, $t_N=T$, with an assumed uniform time step of $\delta = t_{k+1} - t_k$. The extension involves by first introducing two auxiliary intermediate state variables corresponding to the ``{follower}'' and the ``{leader}'' and we then proceed updating the control strategies in a time-parallelized manner.\footnote{Here, we remark that the iterative computations, with appropriate boundary conditions, will involve time-parallelized realization based on multiple shooting methods (see also \cite[Chapter~18]{r8} for additional discussions.)}

\subsubsection{Intermediate state variable w.r.t. the follower} \label{S3.2.1}
Assume that, at the $n^{\rm th}$ iteration, using the admissible control pairs $\bigl(u_1^{(n-1)}(t), u_2^{(n-1)}(t)\bigr)$, we solved both the forward and backward equations w.r.t. the system dynamics of the ``{follower}'' and we further obtained the corresponding solutions at each time instant $t_0=0$, $t_1=\delta$, $t_2=2\delta$, \ldots, $t_N=T$, with uniform time step of $\delta = t_{k+1} - t_k$.

Then, define an intermediate state variable $m_2^{u_2}$ w.r.t. the ``{follower}'' as follows
\begin{align}
m_2^{u_2}(t_k) = \frac{T-t_k}{T} \theta^{u_2}(t_k) + \frac{t_k}{T} p_2(t_k), \label{Eq3.7}
\end{align}
where $\theta^{u_2}(0) = \theta_0$ and $p_2(T) = 0$. Here, we can introduce a family of auxiliary optimal control problems in each of the sub-intervals $[t_k, t_{k+1}]$, with $k \in \{0,1,\ldots,N-1\}$, and we then proceed with minimization of the following 
\begin{align*}
J_2^{k}\bigl [u_2^{(n)} \vert m_2^{u_2} \bigr] &= \frac{1}{2} \bigl \Vert m_2^{u_2}(t_{k+1}) - \theta^{u_2}(t_{k+1}) \bigr \Vert^2 \\
                                                                        & \quad  \quad + \int_{t_k}^{t_{k+1}} \left \{ \frac{\bar{\alpha}}{2} \bigl \Vert \theta^{u_2}(t) \bigr \Vert^2 + \frac{\bar{\beta}}{2} \bigl \Vert u_2^{(n)}(t)  \bigr \Vert^2 \right \} dt \quad \to \quad \min_{u_2^{(n)}(t)} \\
& \quad \quad  \\
 \text{s.t.} \quad &\\
& \quad \quad  \dot{\theta}^{u_2}(t) =-\nabla J_0(\theta^{u_2}(t), \mathcal{Z}^{(1)}) + u_1^{(n-1)}(t_k) + u_2^{(n)}(t), 
\end{align*}
starting from an initial $\theta^{u_2^{(n)}}(t_k)=m_2^{u_2^{(n)}}(t_k)$, while the parameters $\bar{\alpha}$ and $\bar{\beta}$ are given by the following relations $\bar{\alpha} = \frac{\delta}{T} \alpha$ and $\bar{\beta} = \frac{\delta}{T} \beta$, respectively.

Note that at this stage we can establish the relationship between the above decomposed optimization and the original optimization problem. Given an arbitrary intermediate trajectory $m_2^{u_2}(t)$, then we can define the total cost functional on the time interval $[0, T]$ as follows
\begin{align}
\bar{J}_2\bigl [u_2 \vert m_2^{u_2} \bigr] = \frac{T}{\delta} \sum\nolimits_{k=0}^{N-1} J_2^{k}\bigl [u_2^{(n)} \vert m_2^{u_2} \bigr]. \label{Eq3.8}
\end{align}
Then, for arbitrary control strategy pairs $\bigl(u_1(t), u_2(t)\bigr)$, we have 
\begin{align}
m_2^{u_2}(t) \in \argmin_{m_2(t)} \bar{J}_2\bigl [u_2 \vert m_2^{u_2} \bigr]. \label{Eq3.9}
\end{align}
Moreover, the following conditions hold true
\begin{align*}
\bar{J}_2\bigl [u_2 \vert m_2^{u_2} \bigr] = J_2\bigl [u_2\bigr] \quad \text{and} \quad \nabla J_2\bigl [u_2\bigr]\bigl \vert_{[t_k, t_{k+1}]} = \frac{T}{\delta} \nabla J_2^{n}\bigl [u_2 \vert_{[t_k, t_{k+1}]} \vert m_2^{u_2} \bigr].
\end{align*}
Furthermore, the last two equations provide an interpretation for time-parallelized distributed computational implementation, since the intermediate state variable $m_2^{u_2}(t_k)$, with $k=0$, $1$, $2$, \ldots, $N-1$, will allow us to compute independently the corresponding gradient on each time subinterval $[t_k, t_{k+1}]$. Note that part of the modification w.r.t. the ``{follower}''  in the nested algorithm will necessitate the following steps:

\begin{itemize}
\item[{ i.}] Using the admissible control pairs $\bigl(u_1^{(n-1)}(t), u_2^{(n-1)}(t)\bigr)$, solve the forward and backward-equations w.r.t. the system dynamics of the ``{follower}'' and obtain the solutions at time instants $t_0=0$, $t_1=\delta$, $t_2=2\delta$, \ldots, $t_N=T$, with uniform time step of $\delta = t_{k+1} - t_k$.
\item[{ ii.}] Compute the intermediate state variable $m_2^{u_2}$ w.r.t. the ``{follower}'' using
\begin{align*}
m_2^{u_2}(t_k) = \frac{T-t_k}{T} \theta^{u_2}(t_k) + \frac{t_k}{T} p_2(t_k),
\end{align*}
for the time instants $t_0=0$, $t_1=\delta$, $t_2=2\delta$, \ldots, $t_N=T$, with uniform time step of $\delta = t_{k+1} - t_k$.
\item[{ iii.}] On each time subinterval $[t_k, t_{k+1}]$, with $k=0$, $1$, \ldots, $N-1$, determine $u_2^{(n)}(t) \vert_{t \in [t_k, t_{k+1}]}$ as an approximate optimal solution to the problem $\min J_2^{k}\bigl [u_2^{(n)}]$.
\item[{ iv.}] Then, concatenate the controls $u_2^{(n)}(t)\vert_{t \in [t_k, t_{k+1}]}$, for $k=0$, $1$, \ldots, $N-1$, to form $u_2^{(n)}(t)$ for $t \in [0, T]$.
\end{itemize}

\subsubsection{Intermediate state variable w.r.t. the leader} \label{S3.2.2}
Here, we also assume that, at the $n^{\rm th}$ iteration, using the admissible control pairs $\bigl(u_1^{(n-1)}(t), u_2^{(n)}(t)\bigr)$, we have solved both the forward and the backward equations w.r.t. the system dynamics of the ``{leader}'' and we further obtained the corresponding solutions at time instants $t_0=0$, $t_1=\delta$, $t_2=2\delta$, \ldots, $t_N=T$, with uniform time step of $\delta = t_{k+1} - t_k$.

Then, we define an intermediate state variable $m_1^{u_1}$ w.r.t. the ``{leader}'' as follows
\begin{align}
m_1^{u_1}(t_k) = \frac{T-t_k}{T} \theta^{u_1}(t_k) + \frac{t_k}{T} p_1(t_k), \label{Eq3.10}
\end{align}
where $\theta^{u_1}(0) = \theta_0$ and $p_1(T) = -\frac{\partial \Phi(\theta,\mathcal{Z}^{(1)})} {\partial \theta} \biggr \vert_{\theta=\theta^{u_1}(T)}$

Next, we introduce a family of auxiliary optimal control problems in each of the sub-intervals $[t_k, t_{k+1}]$, with $k \in \{0,1,\ldots,N-1\}$, and we then proceed with minimizing 
\begin{align*}
J_1^{k}\bigl [u_1^{(n)} \vert m_1^{u_1} \bigr] &=  \frac{\delta}{2T} \int_{t_k}^{t_{k+1}}  \bigl \Vert \theta^{u_1}(t) \bigr \Vert^2 dt \quad \to \quad \min_{u_1^{(n)}(t)} \\
& \quad \text{s.t.} \quad m_1^{u_1}(t_{k+1}) = -\frac{\partial \Phi(\theta,\mathcal{Z}^{(1)})} {\partial \theta} \biggr \vert_{\theta=\theta^{u_1}(t_{k+1})}&\\
& \quad \quad  \dot{\theta}^{u_1}(t) =-\nabla J_0(\theta^{u_1}(t), \mathcal{Z}^{(1)}) + u_1^{(n)}(t) + u_2^{(n)}(t_k), 
\end{align*}
starting from $\theta^{u_1^{(n)}}(t_k)=m_1^{u_1^{(n)}}(t_k)$. Then, part of the modification w.r.t. the ``{leader}'' in the nested algorithm will necessitate the following steps:

\begin{itemize}
\item[{ i.}] Using the admissible control pairs $\bigl(u_1^{(n-1)}(t), u_2^{(n)}(t)\bigr)$, solve the forward and backward-equations w.r.t. the system dynamics of the ``{follower}'' and obtain the solutions at time instants $t_0=0$, $t_1=\delta$, $t_2=2\delta$, \ldots, $t_N=T$, with uniform time step of $\delta = t_{k+1} - t_k$.
\item[{ ii.}] Compute the intermediate state variable $m_1^{u_1}$ w.r.t. the ``{leader}'' using
\begin{align*}
m_1^{u_1}(t_k) = \frac{T-t_k}{T} \theta^{u_1}(t_k) + \frac{t_k}{T} p_1(t_k),
\end{align*}
for the time instants $t_0=0$, $t_1=\delta$, $t_2=2\delta$, \ldots, $t_N=T$, with uniform time step of $\delta = t_{k+1} - t_k$.
\item[{ iii.}] On each time subinterval $[t_k, t_{k+1}]$, with $k=0$, $1$, \ldots, $N-1$, determine $u_1^{(n)}(t) \vert_{t \in [t_k, t_{k+1}]}$ as an approximate optimal solution to the problem $\min J_1^{k}\bigl [u_1^{(n)}]$.
\item[{ iv.}] Then, concatenate the controls $u_1^{(n)}(t)\vert_{t \in [t_k, t_{k+1}]}$, for $k=0$, $1$, \ldots, $N-1$, to form $u_1^{(n)}(t)$ for $t \in [0, T]$.
\end{itemize}

\subsubsection*{Algorithm -- Intermediate state methods}
Note that if we update the original nested algorithm based on the parts of the above two modifications w.r.t. the ``{follower}'' and the ``{leader}''. Then, we have the following modified nested algorithm based on the intermediate state methods that lends the admissible control updating steps to be fully time-parallelized and, as a result, it provides better the computationally efficiency.

{\rm \footnotesize

{\bf ALGORITHM - 2:} Modified Nested Algorithm Based on the Intermediate State Methods
\begin{itemize}
\item[{\bf 0.}] {\bf Initialize:} Start with any admissible control strategy $u_1^{(n-1)}(t) \in \mathcal{U}_1$ (i.e., an initial guess control function on $[0,T]$) for the ``{\it leader}.'' 
\item[{\bf 1.}] {\bf The regularization-type problem:} Choose an admissible control strategy $u_2^{(n-1)}(t) \in \mathcal{U}_2$ (i.e., an initial guess control function on $[0,T]$) for the ``{\it follower}.'' 
\item[{\bf 2.}] Using the admissible control pairs $\bigl(u_1^{(n-1)}(t), u_2^{(n-1)}(t)\bigr)$, solve the forward and backward-equations w.r.t. the system dynamics of the ``{\it follower},'' i.e., Equations~\eqref{Eq2.9} and \eqref{Eq2.10}, and obtain the solutions at time instants $t_0=0$, $t_1=\delta$, $t_2=2\delta$, \ldots, $t_N=T$, with uniform time step of $\delta = t_{k+1} - t_k$.
\item[{\bf 3.}] Compute the intermediate state variable $m_2^{u_2}$ w.r.t. the ``{follower}'' using
\begin{align*}
m_2^{u_2}(t_k) = \frac{T-t_k}{T} \theta^{u_2}(t_k) + \frac{t_k}{T} p_2(t_k),
\end{align*}
for the time instants $t_0=0$, $t_1=\delta$, $t_2=2\delta$, \ldots, $t_N=T$, with uniform time step of $\delta = t_{k+1} - t_k$.
\item[{\bf 4.}] On each time subinterval $[t_k, t_{k+1}]$, with $k=0$, $1$, \ldots, $N-1$, determine $u_2^{(n)}(t) \vert_{t \in [t_k, t_{k+1}]}$, in time-parallelized manner, as an approximate optimal solution to the problem $\min J_2^{k}\bigl [u_2^{(n)}(\cdot)]$.
\item[{\bf 5.}] Then, concatenate the controls $u_2^{(n)}(t)\vert_{t \in [t_k, t_{k+1}]}$, for $k=0$, $1$, \ldots, $N-1$, to form $u_2^{(n)}(t)$ for $t \in [0, T]$.
 \item[{\bf 6.}] {\bf The controllability-type problem:} Using the admissible control strategy pairs $\bigl(u_1^{(n-1)}(t), u_2^{(n)}(t)\bigr)$, with the updated control strategy $u_2^{(n)}(t)$, solve the forward and backward-equations w.r.t. the system dynamics of the ``{\it leader,}'' i.e., Equations~\eqref{Eq2.17} and \eqref{Eq2.18}, and obtain the solutions at time instants $t_0=0$, $t_1=\delta$, $t_2=2\delta$, \ldots, $t_N=T$, with uniform time step of $\delta = t_{k+1} - t_k$.
 \item[{\bf 7.}] Compute the intermediate state variable $m_1^{u_1}$ w.r.t. the ``{leader}'' using
\begin{align*}
m_1^{u_1}(t_k) = \frac{T-t_k}{T} \theta^{u_1}(t_k) + \frac{t_k}{T} p_1(t_k),
\end{align*}
for the time instants $t_0=0$, $t_1=\delta$, $t_2=2\delta$, \ldots, $t_N=T$, with uniform time step of $\delta = t_{k+1} - t_k$.
\item[{\bf 8.}] On each time subinterval $[t_k, t_{k+1}]$, with $k=0$, $1$, \ldots, $N-1$, determine $u_1^{(n)}(t) \vert_{t \in [t_k, t_{k+1}]}$, in time-parallelized manner, as an approximate optimal solution to the problem $\min J_1^{k}\bigl [u_1^{(n)}(\cdot)]$.
\item[{\bf 9.}] Then, concatenate the controls $u_1^{(n)}(t)\vert_{t \in [t_k, t_{k+1}]}$, for $k=0$, $1$, \ldots, $N-1$, to form $u_1^{(n)}(t)$ for $t \in [0, T]$.
\item[{\bf 10.}] With the updated admissible control strategy pairs $\bigl(u_1^{(n)}(t), u_2^{(n)}(t)\bigr)$, repeat Steps $2$ through $9$, until convergence, i.e., $(1/T)\sum\nolimits_{k=0}^{N-1} \int_{t_k}^{t_{k+1}} \bigl\Vert \nabla J_1^{k}\bigl [u_1^{(n)}(t)]\vert_{t \in [t_k, t_{k+1}]} \bigr\Vert dt \le \epsilon_{\rm tol}$, for some error tolerance $\epsilon_{\rm tol} > 0$.
 \item[{\bf 11.}] {\bf Output:} Return the optimal estimated parameter value $\theta^{\ast} = \theta^{u^{(n)}}(T)$, with $u^{(n)}(t) = (u_2^{(n)}(t), u_2^{(n)}(t))$.
\end{itemize}}

 \begin{remark} \label{R3}
 \begin{itemize}
\item[{ i.}] The computational tasks in Steps~$3$ and $4$ (as well as in Steps~$7$ and $8)$ can be fully time-time parallelized, i.e., the computational tasks can realized in a distributed computational environment with multi-core processors, with possible coarse and fine numerical simulation.  
\item[{ ii.}] In Step~$10$ of the above algorithm, the computational progress can be checked by computing the total average gradients $(1/T)\sum\nolimits_{k=0}^{N-1} \int_{t_k}^{t_{k+1}} \bigl\Vert \nabla J_1^{k}\bigl [u_1^{(n)}(t)]\vert_{t \in [t_k, t_{k+1}]} \bigr\Vert$, which is expected to approach to zero, when $\bigl(u_1^{(n)}(\cdot), u_2^{(n)}(\cdot)\bigr) \to \bigl(u_1^{\ast}(\cdot), u_2^{\ast}(\cdot)\bigr)$, with $u_2^{\ast}(\cdot) \mathcal{F}[u_1^{\ast}(\cdot)]$, as $n \to \infty$ (cf. Remark~\ref{R1}).
\end{itemize}
 \end{remark}
 
\section*{Appendix: Nested Algorithm based-on Successive Approximation Methods}
Here, for the sake of completeness, we provide the nested algorithm, based on successive approximation methods, for numerically solving the optimality conditions of the hierarchical optimal control problem (see \cite{r1} for additional discussions).

{\rm \footnotesize

{\bf ALGORITHM - O:} Nested Algorithm (with Hierarchical Structure) Based-on Successive Approximation Methods
\begin{itemize}
\item[{\bf 0.}] {\bf Initialize:} Start with any admissible control strategy $u_1^{(n-1)}(t) \in \mathcal{U}_1$ (i.e., an initial guess control function on $[0,T]$) for the ``{\it leader}.'' 
\item[{\bf 1.}] {\bf The regularization-type problem:} Choose an admissible control strategy $u_2^{(n-1)}(t) \in \mathcal{U}_2$ (i.e., an initial guess control function on $[0,T]$) for the ``{\it follower}.'' 
\item[{\bf 2.}] Then, using the admissible control pairs $\bigl(u_1^{(n-1)}(t), u_2^{(n-1)}(t)\bigr)$, solve the forward and backward-equations w.r.t. the system dynamics of the ``{\it follower},'' i.e., Equations~\eqref{Eq2.9} and \eqref{Eq2.10}.
\item[{\bf 3.}] Compute the correction term $\delta u_2(\cdot)$ on $[0,T]$ w.r.t. the admissible control $u_2$ of the ``{\it follower}'' (cf. Equation~\eqref{Eq2.11}) using 
 \begin{align*}
  \delta u_2(\cdot) = \gamma_2 \frac{\partial H_2}{\partial u_2} \quad {\text on} \quad [0,T],
\end{align*}
where $\gamma_2 \in (0, 1]$.
 \item[{\bf 4.}] Update the new admissible control for the ``{\it follower}, i.e., $u_2^{(n)}(t)$, using
\begin{align*}
 u_2^{(n)} (t) = u_2^{(n-1)} (t) + \delta u_2(t).
\end{align*}
 \item[{\bf 5.}] {\bf The controllability-type problem:} Using the admissible control strategy pairs $\bigl(u_1^{(n-1)}(t), u_2^{(n)}(t)\bigr)$, with the updated control strategy $u_2^{(n)}(t)$, solve the forward and backward-equations w.r.t. the system dynamics of the ``{\it leader,}'' i.e., Equations~\eqref{Eq2.17} and \eqref{Eq2.18}.
 \item[{\bf 6.}] Compute the correction term $\delta u_1(\cdot)$ on $[0,T]$ w.r.t. the admissible control $u_1$ of the ``{\it leader}'' (cf. Equation~\eqref{Eq2.19}) using 
\begin{align*}
\delta u_1(\cdot) = \gamma_1 \frac{\partial H_1}{\partial u_1} \quad {\text on} \quad [0,T], 
\end{align*}
where $\gamma_1 \in (0, 1]$.
 \item[{\bf 7.}] Update the new admissible control for the ``{\it leader},''  i.e., $u_1^{(n)}(t)$, using
 \begin{align*}
 u_1^{(n)}(t) = u_1^{(n-1)}(t) + \delta u_1(t).
 \end{align*}
\item[{\bf 8.}] With the updated admissible control strategy pairs $\bigl(u_1^{(n)}(t), u_2^{(n)}(t)\bigr)$, repeat Steps $2$ through $7$, until convergence, i.e., $\bigl\Vert \partial H_1/ \partial u_1 \bigr\Vert \le \epsilon_{\rm tol}$, for some error tolerance $\epsilon_{\rm tol} > 0$.
 \item[{\bf 10.}] {\bf Output:} Return the optimal estimated parameter value $\theta^{\ast} = \theta^{u^{(n)}}(T)$, with $u^{(n)}(t) = (u_2^{(n)}(t), u_2^{(n)}(t))$.
\end{itemize}}

\end{document}